\documentclass{birkmult}

\usepackage{amsmath,amsthm,amsfonts,amscd,latexsym,amssymb,amsbsy,dsfont,mathrsfs}
\newtheorem{theorem}{Theorem}[section]
\newtheorem{theorem*}{Theorem}

\newtheorem{proposition}[theorem]{Proposition}

\theoremstyle{definition}
\newtheorem{definition}[theorem]{Definition}
\newtheorem{example}[theorem]{Example}

\theoremstyle{remark}
\newtheorem{remark}[theorem]{Remark}

\numberwithin{equation}{section}

\newcommand{\C}{\mbb{C}}
\newcommand{\OO}{\Omega}

\newcommand{\LL}{\mc{L}}

\newcommand{\K}{\mc{K}}
\newcommand{\mbb}{\mathbb}
\newcommand{\mc}{\mathcal}
\newcommand{\mi}{\mathit}

\newcommand{\mr}{\mathrm}

\newcommand{\lra}{\longrightarrow}

\newcommand\ssp{\sigma_{S}}	%spherical spectrum
\newcommand\srho{\rho_{S}}	%spherical resolvent set

\newcommand{\cS}{\mbb{S}}

\newcommand{\rr}{|}

\newcommand{\1}{\mr{I}}

\newcommand{\cN}{\mc{N}}

\newcommand{\sH}{\mathsf{H}}
\newcommand{\bC}{{\mathbb C}}           %%%  complex numbers
\newcommand{\bH}{{\mathbb H}}
\newcommand{\bN}{{\mathbb N}}
\newcommand{\bR}{{\mathbb R}}
\newcommand{\bS}{{\mathbb S}}

\newcommand{\gB}{{\mathfrak B}}

\def\beq{\begin{equation}}
\def\eeq{\end{equation}}

\def\b{\langle}
\def\k{\rangle}

\begin{document}

%\title[On  compact  normal quaternionic operators]{On  compact quaternionic operators}
\title{Spectral properties of compact normal quaternionic operators}

%%%%%%%

\author{Riccardo Ghiloni}
\address{Department of Mathematics, University of Trento, I--38123, Povo-Trento, Italy}
%\curraddr{}
\email{ghiloni@science.unitn.it}
\thanks{Work partially supported by GNSAGA and  GNFM of INdAM}

\author{Valter Moretti}
\address{Department of Mathematics, University of Trento, I--38123, Povo-Trento, Italy}
%\curraddr{}
\email{moretti@science.unitn.it}
%\thanks{Work partially supported by GNSAGA and  GNFM of INdAM}

\author{Alessandro Perotti}
\address{Department of Mathematics, University of Trento, I--38123, Povo-Trento, Italy}
%\curraddr{}
\email{perotti@science.unitn.it}
%\thanks{Work partially supported by GNSAGA and  GNFM of INdAM}

%%%%%%%


\date{}

%%%%%%%

\begin{abstract}
General, especially spectral, features of compact normal operators in quaternionic Hilbert spaces are studied and some results are established which generalize well-known properties of compact normal operators in complex Hilbert spaces. More precisely, it is proved that 
the norm of such an operator always coincides with the maximum of the set of absolute values of the eigenvalues (exploiting the notion of spherical eigenvalue). Moreover the structure of the spectral decomposition of a generic compact normal operator $T$ is discussed also proving a spectral characterization theorem for compact normal operators.
\end{abstract}

\maketitle

%%%%%%%

%\setcounter{tocdepth}{1}
%\tableofcontents

%%%%%%%
\section{Introduction}

\medskip
Theory of linear operators in quaternionic Hilbert spaces is a well established topic of functional analysis with many applications in physics, especially quantum mechanics (see the introduction of \cite{GhMoPe} for a wide discussion). 
As in complex functional analysis, compact operators play a relevant role as they share features both with generic operators in  infinite dimensional spaces and with matrices in finite dimensional spaces. This intermediate role is particularly evident regarding spectral analysis of normal compact operators. In fact, these operators in infinite dimensional (complex or quaternionic) Hilbert spaces,  on the one hand 
admit a pure point spectrum (except, perhaps, for $0$), on the other hand their spectral expansion needs a proper infinite Hilbertian basis.
This paper is devoted to focus on these peculiar properties exploiting the general framework established in \cite{GhMoPe}.

The notion of spectrum of an operator on quaternionic Hilbert spaces has been introduced only few years ago \cite{libroverde} in the more general context of quaternionic Banach modules. It has been a starting point for developing  functional calculus for the classes of slice and slice regular functions on a quaternionic space (see \cite{libroverde,GhMoPe}).

Let $\sH$ be  a (right) quaternionic Hilbert space (we refer to Section~\ref{HilbertSpaces} for basic defi\-nitions and to \cite{GhMoPe} for more details), let $\gB(\sH)$ be the set of right linear operators on  $\sH$ and  let $\gB_0(\sH)$ be the set of right linear compact operators on  $\sH$.  In \cite{Fashandi2} some properties of compact operators on quaternionic Hilbert spaces were studied. In particular, it was shown that the \emph{spherical spectrum} (cf. Section~\ref{spectrum} for complete definitions) of a compact operator $T$ contains only the eigenvalues of $T$ and possibly 0, and the set of eigenvectors relative to a non-zero eigenvalue $q$ is finite-dimensional.
Another result, similar to what occurs for compact operators in complex Hilbert spaces,  is conjectured in \cite{Fashandi2}:
\paragraph{Conjecture} \textit{If $T\in\gB_0(\sH)$ is self-adjoint, then either $\|T\|$ or $-\|T\|$ is an eigenvalue of $T$.}
\\
In the following we will prove the conjecture for the more general class of normal compact operators on a quaternionic Hilbert space. We also prove the spectral decomposition theorem for normal compact operators and its converse. The complex Hilbert space versions of these results can be found for example in \cite[\S3.3]{Analysisnow}.

\subsection{Main theorems}
As recalled in Section~\ref{spectrum}, the set of eigenvalues of a linear operator $T$ coincides with the \emph{spherical point spectrum}, denoted by $\sigma_{pS}(T)$. We can then rephrase the conjecture in the following way.

\begin{theorem}\label{MT}
Given any normal operator $T \in \gB_0(\sH)$ with spherical point spectrum $\sigma_{pS}(T)$, there exists $\lambda \in \sigma_{pS}(T)$ 
such that:
\beq |\lambda| = \max\{|\mu| \:|\:   \mu \in \sigma_S(T) \}    = \|T\|\:.\label{mi}\eeq
\end{theorem}

The next result is the spectral decomposition for normal compact operators. In the finite-dimensional case, where the compactness requirement is empty, the result is well-known \cite{Jacobson39} (see \cite{FarenickPidkowich} for an ample exposition and further references).

\begin{theorem}\label{MT2}
Given a normal operator $T \in \gB_0(\sH)$ with spherical point spectrum $\sigma_{pS}(T)$, there exists  a Hilbert basis $\cN \subset \sH$ made of eigenvectors of $T$ such that: 
\beq Tx = \sum_{z\in \cN}  z \lambda_z \langle z| x\rangle\quad \mbox{for each $x\in \sH$,} \eeq
where $\lambda_z\in \bH$ is an eigenvalue relative to the eigenvector $z$ and, if $\lambda_z \neq 0$ only a finite number of distinct other  elements $z' \in \cN$ verify 
$\lambda_z = \lambda_{z'}$, moreover the values $\lambda_z$ are at most countably many.

The set $\Lambda$ of eigenvalues $\lambda_z$ with $z\in \cN$ has the property that 
for every $\epsilon>0$ there is a finite set $\Lambda_\epsilon \subset \Lambda$ with
$|\lambda| < \epsilon$ if $\lambda \not \in \Lambda_\epsilon$ (following \cite{Analysisnow} we say that the eigenvalues ``vanish at infinity''). Thus $0$ is the only possible accumulation point of $\Lambda$. If $\sH$ is infinite-dimensional, then $0$ belongs to $\ssp(T)$.
\end{theorem}

\begin{remark} 
Let $\bS$ denote the two-dimensional sphere of imaginary units in $\bH$:
\[\bS:=\{q\in \bH\,|\, q^2=-1\}.\]
For $\imath\in\bS$, let $\bC_\imath$ be the real subalgebra of $\bH$ generated by $\imath$.
We will see in Section \ref{HilbertSpaces} that  to every normal operator $T$ can be associated an anti self-adjoint operator $J$. Along with a chosen imaginary unit $\imath\in\bS$, $J$ defines (cf. Definition \ref{defHJ}) two $\bC_\imath$-Hilbert subspaces of $\sH$, denoted by $\sH^{J\imath}_+$ and  $\sH^{J\imath}_-$.

In the preceding theorem, for every imaginary unit $\imath$, it is possible to choose $\cN$ such that: \beq\{\lambda_z \:| \: z \in \cN\}\setminus \{0\} = (\sigma_{pS}(T)\setminus \{0\}) \cap \bC_\imath\:. \label{spettcond}\eeq
as well as $\cN \subset \sH^{J\imath}_+$.
With a conjugation, in general losing the condition $\cN \subset \sH^{J\imath}_+$,
one can always have $\lambda_z \in \bC_\imath^+$ but $\cN \subset  \sH^{J\imath}_+ \cup  \sH^{J\imath}_-$.
%(2)\label{rem2}  
%(Si dimostra passando a $\sH^{J\imath}_+$  ed usando la stessa propriet\`a per gli operatori compatti su spazi complessi \cite{Analysisnow}, questa in realt\`a \`e una propriet\`a dello spettro degli operatori compatti quaternionici che forse Fashandi non dice)
\end{remark}

%Referring to Remark~\ref{rem2}(2) also the following result holds.
The spectral theorem for compact operators has the following converse. To state it, we need to recall a definition. Given a subset $\K$ of $\bC$, we define the \emph{circularization $\OO_\K$ of $\K$ $($in $\bH \, )$} by setting
\beq \label{circularization}
\OO_\K:=\{\alpha+\jmath\beta \in \bH \,|\, \alpha,\beta \in \bR, \alpha+i\beta \in \K, \jmath \in \bS\}.
\eeq

\begin{theorem}\label{MT3}
Let $T \in \gB(\sH)$. Assume that  there exist a Hilbert basis $\cN$ of $\sH$  and a map $\cN\ni z\mapsto\lambda_z\in\bH$
%an associated set $\Lambda$ of quaternions $\lambda_z$, $z\in \cN$ 
satisfying the following requirements:
\begin{enumerate}
\item[(i)]
$Tx = \sum_{z\in \cN}  z \lambda_z \langle z| x\rangle$ for every $x\in \sH$.
\item[(ii)] For every $z\in\cN$ such that $\lambda_z \neq 0$, only a finite number of distinct other  elements $z' \in \cN$ verify 
$\lambda_z = \lambda_{z'}$;

\item[(iii)] The set $\Lambda$ is countable at most;

\item[(iv)]
For every $\epsilon>0$, there is a finite set $\Lambda_\epsilon \subset \Lambda$ with
$|\lambda| < \epsilon$ if $\lambda \not \in \Lambda_\epsilon$.
\end{enumerate}
Under these conditions $T$ is normal and compact and 
\[  \sigma_{S}(T) \setminus \{0\} = \Omega_\Lambda \setminus \{0\} 
%= \sigma_{pS}(T) \setminus \{0\}
\:.\]
\end{theorem}

\begin{remark}
The structure of the whole spherical spectrum (see Definition \ref{def_spectrum}) of a compact operator $T \in \gB_0(\sH)$ has been studied in \cite[Corollary~2]{Fashandi2}:
\[\sigma_S(T) \setminus\{0\}= \sigma_{pS}(T) \setminus\{0\}\:.\]
If $T$ is normal, then its  \emph{spherical residual spectrum}  (cf.~Section~\ref{spectrum} for definitions) is empty. Therefore in this case if  $0 \in \sigma_{S}(T)\setminus \sigma_{pS}(T)$ then $0$  belongs to the \emph{spherical continuous  spectrum} $\sigma_{cS}(T)$.
\end{remark}

%The last proposition concerns the structure of the whole spherical spectrum (see Definition \ref{def_spectrum}) of a compact normal operator (cf.~also \cite[Corollary~2]{Fashandi2}).
%
%\begin{theorem}\label{MT4}
%Given a normal operator $T \in \gB_0(\sH)$, one has
%\[\sigma_S(T) \setminus\{0\}= \sigma_{pS}(T) \setminus\{0\}\:,\]
%with $0 \in  \sigma_{cS}(T)$ if  $0 \in \sigma_{S}(T)\setminus \sigma_{pS}(T)$.
%\end{theorem}

%\subsection{Summary on mathematical structures} \label{intromath} 

\section{Quaternionic Hilbert spaces}\label{HilbertSpaces}
    We recall some basic notions about quaternionic Hilbert spaces (see e.g.~\cite{BDS}). Let $\bH$ denote the skew field of quaternions.
    Let $\sH$ be a \emph{right} $\bH$--module. $\sH$ is called a \emph{quaternionic pre--Hilbert space} if there exists a Hermitian quaternionic scalar product $\sH\times\sH\ni (u,v)\mapsto \b u|v \k \in \bH$ satisfying the following three properties:
    \begin{itemize}
        \item \emph{Right linearity}: $\b u| vp +  wq \k = \b u | v\k p + \b u | w\k q$ if $p,q \in \bH$ and $u,v,w \in \sH$.
        \item\emph{Quaternionic Hermiticity}: $\b u|v \k = \overline{\b v | u \k}$ if $u,v \in \sH$.
        \item \emph{Positivity}: If $u \in \sH$, then $\b u| u \k \in \bR^+$ and $u=0$ if $\b u | u \k=0$.
    \end{itemize}
    We can define the \emph{quaternionic norm}  by setting
    \[
    \|u\| := \sqrt{\langle u| u\rangle}\in\bR^+ \quad \text{if $u \in \sH$}. 
    \]
    %The function $\|\cdot\|$ is a genuine norm over $\sH$, viewed as a real vector space.
    
    \begin{definition}
        A quaternionic pre--Hilbert space $\sH$ is said to be a \emph{quaternionic Hilbert space} if it is complete with respect to its natural distance $d(u,v):=\|u-v\|$.
    \end{definition}
    
    \begin{example} The space $\bH^n$ with scalar product $\b u,v\k=\sum_{i=1}^n \bar u_i v_i$ is a finite-dimensional quaternionic Hilbert space.\end{example}

Let $\sH$ be a quaternionic Hilbert space.
    \begin{definition}
        A \emph{right $\bH$--linear operator} is a map $T: D(T) \lra \sH$ such that:
        \[
        T(ua+vb)=(Tu)a+(Tv)b \quad \text{if }u,v \in D(T) \text{ and }a,b \in \bH,
        \]
        where the \emph{domain} $D(T)$ of $T$ is a (not necessarily closed) right $\bH$--linear subspace of $\sH$. 
    \end{definition}

       It can be shown that an operator $T: D(T) \lra \sH$ is continuous if and only if it is bounded, i.e.\ there exists $K \geq 0$ such that
    \[
    \|Tu\| \leq K \|u\| \quad \text{for each } u \in D(T).
    \]
    % Every bounded operator is closed.
Let $\|T\|:=\sup_{u \in D(T) \setminus \{0\}} \frac{\|Tu\|}{\|u\|}=\inf\{K \in \bR \:|\: \|Tu\| \leq K \|u\| \ \forall u \in D(T)\}$. 
    The set $\gB(\sH)$ of all bounded operators $T:\sH \lra \sH$ is a complete metric space w.r.t.\ the metric $D(T,S):=\|T-S\|$, 

   Many assertions that are valid in the complex Hilbert spaces case, continue to hold for quaternionic operators. We mention the uniform boundedness principle, the open map theorem, the closed graph theorem, the Riesz representation theorem and the polar decomposition of operators.

As in the complex case, a linear operator $T:\sH\rightarrow\sH$ is called \emph{compact} if it maps bounded sequences to sequences that admit convergent subsequences. 
% The set $\gB_0(\sH)$ of compact operators on  $\sH$ is a subset of $\gB(\sH)$. 
We refer to \cite{Fashandi2} for some properties of compact operator on quaternionic Hilbert spaces. In particular,  $\gB_0(\sH)$ is a closed bilateral ideal of  $\gB(\sH)$ and is closed under adjunction (\cite[Theorem~2]{Fashandi2}).

 \subsection{Left scalar multiplications}

It is possible to equip a (right) quaternionic Hilbert space $\sH$ with a \emph{left} multiplication by quaternions. It is a non--canonical operation relying upon a choice of a preferred Hilbert basis. So, pick out a Hilbert basis $\cN$ of $\sH$ and define the \emph{left scalar multiplication of $\sH$ induced by $\cN$} as the map $\bH \times \sH \ni (q,u) \mapsto qu \in \sH$ given by 
     \[
     q u:=\textstyle\sum_{z \in \cN} z q \b z|u \k \quad \text{if $u \in \sH$ and $q \in \bH$.}
     \]
     For every $q \in \bH$, the map $L_q:u\mapsto qu$ belongs to $\gB(\sH)$. Moreover, the  map $\LL_\cN:\bH \lra \gB(\sH)$, defined by setting $\LL_\cN(q):=L_q$ is a norm--preserving real algebra homomorphism. 
     
     The set $\gB(\sH)$ is always a \emph{real Banach $C^*$--algebra with unity}. It suffices to consider the right scalar multiplication $(Tr)(u)=T(u)r$ for real $r$ and the adjun\-ction $T \mapsto T^*$ as $^*$--involution. By means of a left scalar multiplication, it can be given the richer structure of \emph{quaternionic Banach $C^*$--algebra}.
     
     \begin{theorem}[\cite{GhMoPe}\S3.2]
         Let $\sH$ be a quaternionic Hilbert space equipped with a left scalar multiplication. Then the set $\gB(\sH)$, equipped with the pointwise sum, with the scalar multiplications defined by
         \[(qT)u:=q(Tu)\quad and\quad(Tq)(u):=T(qu),\]
         with the composition as product and with  $T \mapsto T^*$ as $^*$-involution, is a \emph{quaternionic two--sided Banach $C^*$--algebra} \emph{with unity}.
     \end{theorem}
      Observe that the map $\LL_\cN$ gives a $^*$-representation of $\bH$ in $\gB(\sH)$.

\subsection{Imaginary units and complex subspaces}

Consider a quaternionic Hilbert space $\sH$ equipped with a left scalar multiplication $\bH \ni q \mapsto L_q$. For short, we write $L_qu=qu$. For every imaginary unit $\imath \in \bS$, the operator $J:= L_\imath$ is anti self--adjoint and unitary; that is, it holds:
\[
J^*=-J \quad \text{and} \quad J^*J=\1.
\]
%The proof straightforwardly arises from Proposition~\ref{propprod}. We intend to establish 
It holds also the converse statement: if an operator $J \in \gB(\sH)$ is anti--self adjoint and unitary, then $J=L'_\imath$ for some left scalar multiplication of $\sH$ (see \cite[Proposition~3.8]{GhMoPe}).

In the following, we also need a definition known from the litera\-ture \cite{Emch}.

\begin{definition}\label{defHJ}
Let $J \in \gB(\sH)$ be an anti self--adjoint and unitary operator and let $\imath\in \bS$. Let $\bC_\imath$ denote the real subalgebra of $\bH$ generated by $\imath$; that is, $\C_\imath:=\{\alpha+\imath\beta \in \bH \, | \, \alpha,\beta \in \bR\}$. Define the $\bC_\imath$-\emph{complex subspaces $\sH^{J\imath}_+$ and $\sH^{J\imath}_-$ of $\sH$ associated with $J$ and $\imath$} by setting
\[
\sH^{J\imath}_\pm:=\{u \in \sH \, | \, Ju=\pm u \imath\}.
\]
\end{definition}

\begin{remark} \label{rC}
$\sH^{J\imath}_\pm$ are closed subsets of $\sH$, because $u \mapsto Ju$ and $u \mapsto \pm u\imath$ are continuous. However, they are not (right $\bH$--linear) subspaces of $\sH$. Note also that the space $\sH$ admits the direct sum decomposition
\[\sH=\sH^{J\imath}_+\oplus \sH^{J\imath}_-,\]
with projections $\sH\ni x\mapsto P_\pm(x):=\frac12(x\mp Jx\imath)\in\sH^{J\imath}_\pm$.
\end{remark}

\subsection{Resolvent and spectrum}\label{spectrum}

    It is not clear how to extend the definitions of spectrum and resolvent in quaternionic Hilbert spaces. Let us focus on the simpler case of eigenvalues of a bounded right $\bH$--linear operator $T$. Without fixing any left scalar multiplication of $\sH$, the equation determining the eigenvalues reads as follows:
    \[
    Tu=uq.
    \] 
    Here a drawback arises: if $q \in \bH \setminus \bR$ is fixed, the map $u \mapsto uq$ is not right $\bH$--linear. Consequently, the eigenspace of $q$ cannot be a right $\bH$--linear subspace. Indeed, if $\lambda \neq 0$, $u\lambda$ is an eigenvector of $\lambda^{-1} q \lambda$ instead of $q$ itself. As a  second guess, one could decide to deal with quaternionic Hilbert spaces equipped with a left scalar multiplication and require that
    \[
    Tu=qu.
    \]
    Now both sides are right $\bH$--linear. However, this approach is not suitable for physical applications, where self--adjoint operators should have real spectrum.    
   We come back to the former approach and accept that each eigenvalue $q$ brings a whole conjugation class of the quaternions, the \emph{eigensphere} \[\cS_q:=\{\lambda^{-1} q \lambda \in \bH \,|\, \lambda \in \bH \setminus \{0\}\}.\]
%    The eigenspace of $\cS_q$ \emph{is} a right $\bH$-linear subspace.
      
    We adopt the viewpoint introduced in \cite{libroverde} for quaternionic two-sided Banach modules.
Given an operator $T:D(T) \lra \sH$ and $q \in \bH$, let \[ \Delta_q(T):=T^2-T(q+\overline{q})+\1 |q|^2.\]
    
    \begin{definition}\label{def_spectrum}
        The \emph{spherical resolvent set} of $T$ is the set $\srho(T) $ of $q\in\bH$ such that:
        \begin{itemize}
            \item[({a})] $\mi{Ker}(\Delta_q(T))= \{0\}$.  
            \item[({b})] $\mi{Range}(\Delta_q(T))$ is dense in $\sH$. 
            \item[({c})] $\Delta_q(T)^{-1}:\mi{Range}(\Delta_q(T)) \lra D(T^2)$ is bounded.
        \end{itemize}  
            \end{definition}
        The \emph{spherical spectrum} $\ssp(T)$ of $T$ is defined by $\ssp(T):=\bH \setminus \rho_S(T)$.     
It decomposes into three disjoint \emph{circular} (i.e.\ invariant by conjugation) subsets:  
    \begin{itemize}
        \item [(i)] the \emph{spherical point spectrum} of $T$ (the set of \emph{eigenvalues}):
        \[
        \sigma_{\mi{pS}}(T):=\{q \in \bH \,|\, \mi{Ker}(\Delta_q(T)) \neq \{0\}\}.
        \]
        \item [(ii)] the \emph{spherical residual spectrum} of $T$:
        \[
        \sigma_{\mi{rS}}(T):=\left\{q \in \bH \left| \, \mi{Ker}(\Delta_q(T))=\{0\}, \, \overline{\mi{Range}(\Delta_q(T))} \neq \sH \right.\right\}.
        \] 
        \item [(iii)] the \emph{spherical continuous spectrum} of $T$:
        \[
        \sigma_{\mi{cS}}(T):=\left\{q \in \bH \left| \,   \Delta_q(T)^{-1} \text{ is densely defined but not bounded} \right.\right\}.
        \]
    \end{itemize}         
The \emph{spherical spectral radius} of $T$ is defined as
\[
r_S(T):=\sup\big\{|q|  \big| \, q \in  \ssp(T)\big\}\in \bR^+ \cup \{+\infty\}.
\]
    In our context, the subspace $\mi{Ker}(\Delta_q(T))$ has the role of an ``eigenspace''. In particular, $\mi{Ker}(\Delta_q(T)) \neq \{0\}$ if and only if $\cS_q$ is an eigensphere of $T$.    

\subsection{Spectral properties}
    
    The spherical resolvent and the spherical spectrum can be defined for bounded right $\bH$--linear operators on quaternionic two--sided Banach modules in a form similar to that introduced above (see \cite{libroverde}). Several spectral properties of bounded operators on complex Banach or Hilbert  spaces remain valid in that general context. Here we recall some of these properties in the quaternionic Hilbert setting (cf.~Theorem~4.3 in \cite{GhMoPe}).
    
    \begin{theorem}[\cite{GhMoPe}\S4.1] \label{teopropspectrum}
        Let $\sH$ be a quaternionic Hilbert space and let $T\in \gB(\sH)$. Then
        \begin{itemize}
            \item[(a)]
            $r_S(T) \leq \|T\|$.
            \item[(b)] $\ssp(T)$ is a non--empty compact subset of $\bH$.
             \item[(c)] Let $P \in \bR[X]$. Then, if $T$ is self--adjoint, the following \emph{spectral map property} holds:
             \[
             \ssp(P(T))=P(\ssp(T)).
             \]
            \item[(d)] \emph{Gelfand's spectral radius formula} holds:
            \[ \label{GF}
            r_S(T)=\lim_{n \to +\infty} \|T^n\|^{1/n}.
            \]
            In particular, if $T$ is \emph{normal} (i.e.\ $TT^*=T^*T$), then  $ r_S(T)=\|T\|$.
        \end{itemize}
    \end{theorem}

 Regardless different definitions with respect to the complex Hilbert space case, the notions of spherical spectrum and resolvent set enjoy some properties which are quite similar to those for complex Hilbert spaces. Other features, conversely, are proper to the quaternionic Hilbert space case.
 First of all, it turns out that the spherical point spectrum coincides with the set of eigenvalues of $T$.
 
 \begin{proposition}\label{propsigmap}
     Let $\sH$ be a quaternionic Hilbert space and let $T:D(T) \lra \sH$ be an operator. Then $\sigma_{pS}(T)$ coincides with the set of all eigenvalues of $T$.
    \end{proposition}

The subspace $\mi{Ker}(\Delta_q(T))$  has the role of an eigenspace of $T$. In particular, $\mi{Ker}(\Delta_q(T))\ne\{0\}$ if and only if $\bS_q$ is an eigensphere of $T$. 
    
    \begin{theorem}\label{teospectra} 
        Let  $T$ be an operator with dense domain on  a quaternionic Hilbert space $\sH$.
        \begin{itemize}
            \item[(a)]   $\ssp(T)=\ssp(T^*)$. 
            \item[(b)] If $T\in \gB(\sH)$ is normal, then
            \begin{itemize}
                \item[(i)] $\sigma_{\mi{pS}}(T)=\sigma_{\mi{pS}}(T^*)$.
                \item[(ii)] $\sigma_{\mi{rS}}(T)=\sigma_{\mi{rS}}(T^*)=\emptyset$.
                \item[(iii)] $\sigma_{\mi{cS}}(T)=\sigma_{\mi{cS}}(T^*)$.
            \end{itemize}
            \item[(c)]  If  $T$ is self--adjoint, then   $\ssp(T) \subset \bR$ and $\sigma_{\mi{rS}}(T)$ is empty.
            \item[(d)] If $T$ is anti self--adjoint,  then $\ssp(T) \subset \mr{Im}(\bH)$ and $\sigma_{\mi{rS}}(T)$ is empty.
            \item[(e)] If $T \in \gB(\sH)$ is unitary,  then $\ssp(T) \subset \{q \in \bH \,|\, |q|=1\}$.
            \item[(f)] If $T\in \gB(\sH)$ is anti self--adjoint and unitary, then $\ssp(T)=\sigma_{pS}(T) = \bS$.
%            (the sphere of quaternionic imaginary units).
        \end{itemize}
    \end{theorem}    
    
    It can be shown that, differently from operators on complex Hilbert spaces, a  normal operator  $T$ on a quaternionic space  is unitarily equivalent to $T^*$.

\subsection{Compact operators}\label{CompactOperators}
Compact operators have some peculiar spectral properties. Some of them were investigated in \cite{Fashandi2}. In particular, if $T\in\gB_0(\sH)$ and $q\in\sigma_{\mi{pS}}(T)\setminus\{0\}$ is an eigenvalue, then $\mi{Ker}(\Delta_q(T))$ has finite dimension \cite[Theorem 3]{Fashandi2}. Moreover, the spherical spectrum of $T\in\gB_0(\sH)$ consists only of  the eigenvalues of $T$ and (possibly) 0 (cf. \cite[Corollary 2]{Fashandi2}):
\[ \ssp(T)\setminus\{0\}=\sigma_{pS}(T)\setminus\{0\}.\]

\subsection{Slice nature of normal operators}

We recall the ``slice'' character of $\bH$: 
        \begin{itemize}
            \item  $\bH=\bigcup_{\jmath \in \bS}\bC_\jmath$\quad where $\bC_\jmath$ is the real subalgebra $\b \jmath\k\simeq\bC$.
            \item  $\bC_\jmath \cap \bC_\kappa=\bR$\quad for every $\jmath,\kappa \in \bS$ with $\jmath \neq \pm\kappa$.
        \end{itemize}
 
This decomposition of $\bH$ has an ``operatorial'' counterpart on a quaternionic Hilbert space.  It was established in Theorems 5.9 of \cite{GhMoPe}.
    
    \begin{theorem}[\cite{GhMoPe}\S5.4]\label{AplusJB}
        Given any normal operator $T \in \gB(\sH)$, there exist three operators $A,B,J \in \gB(\sH)$ such that:
        \begin{itemize}
            \item[({i})] $T=A+JB$.
            \item[({ii})] $A$ is self--adjoint and $B$ is positive.
            \item[({iii})] $J$ is anti self--adjoint and unitary.
            \item[({iv})] $A$, $B$ and $J$ commute mutually.
        \end{itemize}
 
        Furthermore, it holds:
        \begin{itemize}
            \item $A$ and $B$ are uniquely determined by $T$: $A=(T+T^*)\frac{1}{2}$ and $B=|T-T^*|\frac{1}{2}$.
            \item $J$ is uniquely determined by $T$ on $\mi{Ker}(T-T^*)^\perp$.
        \end{itemize}    
    (where for $S\in\gB(\sH)$, $|S|$ denotes the operator defined as the square root of the positive operator $S^*S$).
    \end{theorem}

In the following, we denote by $\sigma(B)$ and $\rho(B)$ the standard spectrum and resolvent set
of a bounded operator $B$ of a complex Hilbert space, respectively.

\begin{proposition}[\cite{GhMoPe}\S5.4]\label{propinterssigma} 
Let $\sH$ be a quaternionic Hilbert space, let $T \in \gB(\sH)$ be a normal operator, let $J \in \gB(\sH)$ be an anti self--adjoint and unitary operator satisfying $TJ=JT$,  $T^*J=JT^*$, let $\imath \in \bS$ and let $\sH_\pm^{J\imath}$ be the complex subspaces of $\sH$ associated with $J$ and $\imath$ (see Definition~\ref{defHJ}). Then we have that
\begin{itemize}
 \item[$(\mr{a})$] $T(\sH_+^{J\imath}) \subset \sH_+^{J\imath}$ and $T^*(\sH_+^{J\imath}) \subset \sH_+^{J\imath}$.
\end{itemize}
Moreover, if $T\rr_{\sH_+^{J\imath}}$ and $T^*\rr_{\sH_+^{J\imath}}$ denote the $\bC_{\imath}$--complex operators in $\gB(\sH_+^{J\imath})$ obtained restricting respectively $T$ and $T^*$ to $\sH_+^{J\imath}$, then it holds:
\begin{itemize}
 \item[$(\mr{b})$] $(T\rr_{\sH_+^{J\imath}})^*= T^*\rr_{\sH_+^{J\imath}}$.
 \item[$(\mr{c})$] $\sigma(T\rr_{\sH_+^{J\imath}}) \cup \overline{\sigma(T\rr_{\sH_+^{J\imath}})}= \ssp(T)\cap \bC_\imath$. Here $\sigma(T\rr_{\sH_+^{J\imath}})$ is considered as a subset of $\bC_\imath$ via the natural identification of $\bC$ with $\bC_\imath$ induced by the real vector isomorphism $\bC \ni \alpha+i\beta \mapsto \alpha+\imath\beta \in \bC_\imath$.
 \item[$(\mr{d})$]
$\ssp(T)=\OO_\K$, where $\K:=\sigma(T\rr_{\sH_+^{J\imath}})$.
\end{itemize}

An analogous statement holds for $\sH^{J\imath}_-$.
\end{proposition}
%%%

\section{Proofs of the main results}

\subsection{Proof of Theorem~\ref{MT}}
If $T =0$ there is nothing to prove, since $0 \in \sigma_{pS}(T)$ and $\|T\|=0$ in that case. So, we henceforth assume that $T\neq 0$. 
Since $T \in \gB(\sH)$ is normal, Theorem~\ref{AplusJB} assures the existence of  an anti self-adjoint unitary 
right $\bH$-linear operator $J : \sH \to \sH$, commuting with $T$ and $T^*$ and fulfilling
$T= (T+T^*)\frac{1}{2} + J|T-T^*|\frac{1}{2}$.  Next, if $\sH_+^{J\imath}$ is the complex subspace 
associated with a imaginary unit $\imath \in \bS$ as in Definition~\ref{defHJ}, it turns out that
$T$ is the unique right $\bH$-linear operator, defined on $\sH$, whose restriction to $\sH_+^{J\imath}$ coincides with the complex-linear operator $S:= T |_{\sH_+^{J\imath}}: \sH_+^{J\imath} \to \sH_+^{J\imath}$ (it immediately arises from
 Propositions~3.11  and 5.11 in \cite{GhMoPe}).
We also know that $\|S\|= \|T\|$, in view of Proposition~3.11 in \cite{GhMoPe}.

By hypotheses $T$ is compact and thus $S$ is compact as well as we go to prove. If $\{u_n\}_{n\in \bN} \subset \sH_+^{J\imath}$ is a bounded sequence of vectors, it is a bounded sequence of vectors of $\sH$ too,  and thus  the sequence $\{Tu_n\}_{n\in \bN}$ admits a subsequence $\{Tu_{n_k}\}_{k\in \bN}$ converging to some   $v \in \sH$, because $T$ is compact.  However, since $\sH_+^{J\imath}$
is closed (because of its definition and the fact that $J$ is continuous), we also have that $v \in \sH_+^{J\imath}$ and  that $\{Su_{n_k}\}_{k\in \bN}$ converges 
to   $v \in \sH_+^{J\imath}$, because $Tu_n = Su_n$. We have found that, 
for every bounded sequence $\{u_n\}_{n\in \bN} \subset \sH_+^{J\imath}$, there is a subsequence of $\{Su_n\}_{n\in \bN}$
converging to some    $v \in \sH_+^{J\imath}$. Thus $S$ is compact.

To go on,  Lemma~3.3.7 in \cite{Analysisnow} entails that there exists $\lambda \in \sigma_p(S)$
with $|\lambda| = \|S\|$. Notice that $\lambda \neq 0$ otherwise $S=0$ and thus $T=0$ %(not permitted)  
by uniqueness of the extension of $S$.
Finally, point (d) of Proposition~\ref{propinterssigma} %  5.11 in \cite{GhMoPe} 
implies that $\lambda \in \sigma_{S}(T)$.
Since $T$ is compact,  by Corollary~2 of \cite{Fashandi2}, we have that
$\lambda \in \sigma_{pS}(T)$. Summing up, we have obtained that there is 
$\lambda \in \sigma_{pS}(T)$ with $|\lambda| = \|S\|= \|T\|$, where the absolute value is, indifferently, that in $\bC$ or that in $\bH$. The remaining identity in (\ref{mi}) is now equivalent to:
$\sup\{|\mu| \:|\:   \mu \in \sigma_S(T) \}    = \|T\|$, i.e.\ $r_S(T)=\|T\|$. In this form, it was proved in point (d) of Theorem~\ref{teopropspectrum}. 

\subsection{Proof of Theorem~\ref{MT2}} Fix an imaginary unit $\imath\in\bS$ and consider the normal compact operator
$S : \sH^{J\imath}_+ \to \sH^{J\imath}_+$ as in the proof of Theorem~\ref{MT}.
As a consequence of Theorem~3.3.8 in \cite{Analysisnow}, there exist a Hilbert basis 
$\cN\subset \sH^{J\imath}_+$ made of eigenvectors of $S$ and a map $\cN\ni z\mapsto\lambda_z\in\bC_\imath$ such that 
 each $\lambda_z$ is an eigenvalue of $S$ in $\bC_\imath$ relative to $z\in\cN$ and, if $\lambda_z \neq 0$, only a finite number of distinct other  elements $z' \in \cN$ verify 
$\lambda_z = \lambda_{z'}$. Moreover, the values $\lambda_z$ are at most countably many.
We know by  Lemma~3.10(b) in \cite{GhMoPe}  that $\cN$ is also a Hilbert basis of $\sH$, so that, 
if $x \in \sH$, then 
\[x = \sum_{z\in\cN} z \langle z | x\rangle \:.\]
Since $T$ is continuous and $Tz= Sz=z \lambda_z$, we have:
 \[Tx = \sum_{z\in \cN} z \lambda_z \langle z | x\rangle \:.\]
>From Theorem~3.3.8 in \cite{Analysisnow}, we also get that the set $\{\lambda_z\}_{z\in\cN}$ of the eigenvalues of $S$, and therefore those of $T$, vanish at infinity.
 Thus $0$ is the only possible accumulation point of $\Lambda$. 
%The proof of the last statement is immediate.
If $\sH$ is not finite-dimensional and $0$ is not an eigenvalue of $T$, then $0$ must be an accumulation point of $\Lambda$, since every set of eigenvectors $\mi{Ker}(\Delta_{\lambda_z}(T))$  has finite dimension if $\lambda_z\ne0$. Since $\ssp(T)$ is closed, in any case $0\in\ssp(T)$.

\begin{remark}
Since $\lambda_z\in\bC_\imath$ for each $z\in\cN$,  equation \eqref{spettcond} follows from Corollary~2 in \cite{Fashandi2}.
\end{remark}

\subsection{Proof of Theorem~\ref{MT3}} 
Fix an imaginary unit $\imath\in\bS$.
For each eigenvector $z\in\cN$, with  non real eigenvalue $\lambda_z\in\bH$, we can choose a unit quaternion $\mu_z$ such that $\mu_z^{-1}\lambda_z\mu_z$ belongs to the intersection of the eigensphere of $\lambda_z$ with the complex plane $\bC_\imath$. This means that $z':=z\mu_z$ is still an eigenvector of $T$, with eigenvalue $\lambda_{z'}:=\mu_z^{-1}\lambda_z\mu_z\in\bC_\imath$. If $\lambda_z\in\bR$, we set $\mu_z=1$. The set $\cN':=\{z\mu_z\, |\, z\in\cN\}$ is still a Hilbert basis of $\sH$, such that $Tx=\sum_{z\in\cN'}z\lambda_z\b z|x\k$ for every $x\in\sH$.

The linear operator $J$ defined by setting
\[Jx:=\sum_{z\in\cN'}z\imath\b z| x\k\]
is an anti self-adjoint and unitary operator on $\sH$ (cf.~Proposition~3.1 in \cite{GhMoPe}).
Since
\[
TJx=\sum_{z\in\cN'}z\lambda_z\b z|Jx\k=\sum_{z\in\cN'}z\lambda_z\imath\b z|x\k\]
and
\[JTx=\sum_{z\in\cN'}z\imath\b z|Tx\k=\sum_{z\in\cN'}z\imath\lambda_z\b z|x\k,
\]
we have that $J$ and $T$ commute. Moreover, since $T^*x=\sum_{z\in\cN'}z\overline{\lambda}_z\b z|x\k$, the same holds for $J$ and $T^*$.  Let $\sH_\pm^{J\imath}$ be the complex subspaces of $\sH$ associated with $J$ and the imaginary unit $\imath \in \bS$ as in Definition~\ref{defHJ}. Observe that $\cN'\subset\sH_+^{J\imath}$, since $Jz=z\imath$ for each $z\in\cN'$. 
Let $\jmath\in\bS$ be an imaginary unit orthogonal to $\imath$. Then $\cN'\jmath:=\{z\jmath\,|\,z\in\cN'\}$ is a Hilbert basis for $\sH_-^{J\imath}$ (cf.~Lemma~3.10 in \cite{GhMoPe}).
The $\bC_\imath$-complex subspaces $\sH_\pm^{J\imath}$ of $\sH$ are invariant for $T$, since
\[JTu=TJu=TJ(\mp Ju\imath)=\pm(Tu)\imath\text{\quad for each }u\in \sH_\pm^{J\imath}.\]
 Let $S_\pm:=T|_{\sH_\pm^{J\imath}}$ be the restrictions of $T$. Then $S_\pm$ are diagonalizable, since \[S_+x_+=\sum_{z\in\cN'}z\lambda_z\b z|x_+\k\text{\quad and\quad}S_-x_-=\sum_{z\in\cN'}z\jmath\overline{\lambda}_z\b z\jmath|x_-\k\]
 for every $x_\pm\in\sH_\pm^{J\imath}$. We can then apply Theorem~3.3.8 of \cite{Analysisnow} and obtain that $S_\pm$ are normal and compact. Using Proposition~3.11 of  \cite{GhMoPe}, we get the normality of $T$. 

It remains to prove that $T$ is compact. 
Recall from Remark~\ref{rC} that $\sH=\sH_+^{J\imath}\oplus\sH_-^{J\imath}$, with projections $P_\pm$ defined by $P_\pm(x)=\frac12(x\mp Jx\imath)\in\sH^{J\imath}_\pm$.
If $\{x_n\}_{n\in \bN} \subset \sH$ is a bounded sequence of vectors, then also $\{P_+x_n\}_{n\in \bN} \subset \sH_+^{J\imath}$ and $\{P_-x_n\}_{n\in \bN} \subset \sH_-^{J\imath}$ are bounded sequences.
Since $S_+$ is compact, the sequence $\{S_+P_+x_n\}_{n\in \bN}$ admits a subsequence $\{S_+P_+x_{n_k}\}_{k\in \bN}$ converging to some  $v_+\in \sH_+^{J\imath}$. Similarly, since $S_-$ is compact, we can extract from the sequence $\{S_-P_-x_{n_k}\}_{k\in \bN}$
a subsequence $\{S_-P_-x_{n_{k_l}}\}_{l\in \bN}$ converging to some  $v_-\in \sH_-^{J\imath}$. Then the sequence $\{Tx_{n_{k_l}}\}_{l\in \bN}$ converges to $v:=v_++v_-$, since  $Tx=S_+P_+x+S_-P_-x$ for each $x\in\sH$.

We have shown that for every bounded sequence $\{x_n\}_{n\in \bN} \subset \sH$ there is a subsequence of $\{Tx_n\}_{n\in \bN}$ converging to some    $v \in \sH$. Thus $T$ is compact.

The last statement concerning the spectrum of $T$ follows from Proposition~\ref{propinterssigma}. %5.11 in \cite{GhMoPe}.

\end{document}